\documentclass[final,1p,a4paper,10pt]{elsarticle}
\usepackage{caption,subcaption}
\usepackage{graphicx}
\usepackage{epstopdf}
\usepackage{color}
\usepackage{amsmath}
\usepackage{amssymb}
\usepackage{units}
\usepackage{float}
\usepackage{natbib}
\usepackage{amsthm}
\usepackage{amsfonts}
\usepackage{epsfig}
\usepackage{mathtools}
\usepackage{booktabs,lipsum}
\usepackage[export]{adjustbox}
\usepackage{verbatim}
\usepackage[labelsep=period]{caption}
\usepackage{multirow}
\usepackage{prettyref}
\newrefformat{tab}{\tablename~\ref{#1}}
\usepackage{lineno}
\usepackage[colorlinks=true,linkcolor=blue]{hyperref}
\biboptions{square,numbers,sort&compress}
\journal{??}
\begin{document}
\begin{frontmatter}
\title{RBF-DQ Method for Solving Non-linear Differential Equations of Lane-Emden type}
\author[mymainaddress]{K. Parand\corref{mycorrespondingauthor}}
\cortext[mycorrespondingauthor]{Corresponding author. Member of research group of Scientific Computing. Fax:
+98 2122431653.}
\ead{k\_parand@sbu.ac.ir}
\author[mymainaddress]{S. Hashemi}
\ead{soleiman.hashemi.sh@gmail.com}
\address[mymainaddress]{Department of Computer Sciences, Shahid Beheshti University, G.C. Tehran 19697-64166, Iran}
\begin{abstract}
The Lane-Emden type equations are employed in the modelling of several phenomena in the areas of mathematical physics and astrophysics. In this paper, a new numerical method is applied to investigate some well-known classes of Lane-Emden type equations which are non-linear ordinary differential equations on the semi-infinite domain $[0, \infty)$. We will apply a meshless method based on radial basis function differential quadrature (RBF-DQ) method. In RBFs-DQ, the derivative value of function with respect to a point is directly approximated by a linear combination of all functional values in the global domain. The main aim of this method is the determination of weight coefficients. Here, we concentrate on Gaussian(GS) as a radial function for approximating the solution of the mentioned equation. The comparison of the results with the other numerical methods shows the efficiency and accuracy of this method. 
\end{abstract}
\end{frontmatter}
\section{Introduction}
\subsection{\textbf{The Lane-Emden equation}}
The study of a singular initial value problems modelled by second-order non-linear ordinary differential equations (ODEs) on a
semi-infinite domain have attracted many mathematicians and physicists. One of the equations in this category is the following Lane-Emden type equation:
\begin{equation}\label{Lane-Emden}
y''(x) + \frac {\alpha}{x} y'(x) + f(x,y(x)) = h(x), \texttt{ $\alpha x \geq 0 $,}
\end{equation}
with initial conditions:
\begin{eqnarray}\label{conditions}
&&y(0) = A,\\\nonumber
&&y'(0)= B,
\end{eqnarray}
where $\alpha$, $A$ and $B$ are real constants and $f(x,y)$ and $h(x)$ are some given continuous real-valued functions.
For special forms of $f (x, y)$, the well-known Lane-Emden equations occur
in several models of non-Newtonian fluid mechanics, mathematical physics, astrophysics, etc \citep{Chandrasekhar,Davis}. For example, $f (x, y) = q(y)$, the Lane-Emden equations occur in modelling several phenomena in mathematical physics and
astrophysics such as the theory of stellar structure, the thermal behaviour of a spherical cloud of gas, isothermal gas sphere and
theory of thermionic currents \citep{Chandrasekhar,Davis}.\\
Let $P(r)$ denote the total pressure at a distance $r$ from the center of spherical gas cloud. The total pressure is due to the usual gas pressure and a contribution from radiation:
\begin{equation*}
P=(\frac{1}{3}\xi T^{4}+\frac{R T}{\upsilon}),
\end{equation*}
where $\xi$, $T$, $R$ and $\upsilon$ are respectively the radiation constant, the absolute temperature, the gas constant, and the specific volume, respectively \citep{Agarwala.DonalORegan,Chandra}. Let $M(r)$ be the mass within a sphere of radius $r$ and $G$ the constant of gravitation. The equilibrium equation for the configuration are
\begin{equation}\label{la2}
\frac{d P}{d r}= - \rho\frac{G M(r)}{r^{2}},
\end{equation}
\begin{equation*}
\frac{d M(r)}{d r}=4 \pi \rho r^{2},
\end{equation*}
where $\rho$, is the density, at a distance $r$, from the center of a spherical star \citep{Chandrasekhar}.
Eliminating $M$ yields of these equations results in the following equation, which as should
be noted, is an equivalent form of the Poisson equation \citep{Chandra, Horedt}:
\begin{equation*}
\frac{1}{r^{2}}\frac{d}{d r}(\frac{r^{2}}{\rho}\frac{d P}{d r})=-4\pi G \rho.
\end{equation*}
We already know that in the case of a degenerate electron gas, the pressure and density are $\rho =P^{3/5}$, assuming that such a relation
exists for other states of the star, we are led to consider a relation of the form $P=K \rho^{1+\frac{1}{m}}$, where $K$ and $m$ are constants.
\par
We can insert this relation into Eq.\eqref{la2} for the hydrostatic equilibrium condition and, from this, we can rewrite the equation as follows:
\begin{equation*}
\big[\frac{K(m+1)}{4\pi G}\lambda^{\frac{1}{m}-1}\big]\frac{1}{r^{2}}\frac{d}{d r}(r^{2}\frac{d y}{d r})=-y^{m},
\end{equation*}
where $\lambda$ represents the central density of the star and $y$ denotes the dimensionless quantity, which are both related to $\rho$ through the following relation \citep{Chandrasekhar, Horedt}:
\begin{equation*}
\rho = \lambda y^{m}(x),
\end{equation*}
and let
\begin{equation*}
r=ax,
\end{equation*}
\begin{equation*}
a=\Bigg[\frac{K(m+1)}{4\pi G}\lambda^{\frac{1}{m}-1}\Bigg]^{\frac{1}{2}}.
\end{equation*}
Inserting these relations into our previous relation we obtain the Lane-Emden equation \citep{Chandra, Horedt}:
\begin{equation*}
\frac{1}{x^{2}}\frac{d}{dx}(x^{2}\frac{d y}{d x})=- y^{m},
\end{equation*}
now, we will have the standard Lane-Emden equation with $f(x,y)=y^m$
\begin{equation}\label{Lane-Emden equation}
y''(x)+\frac{2}{x}y'(x)+{y^{m}(x)}=0 ~~ , ~~~ x>0 ,
\end{equation}
at this point, it is also important to introduce the boundary conditions which are based upon the following boundary conditions
for hydrostatic equilibrium and normalization consideration of the newly introduced quantities $x$ and $y$. What follows for $r = 0$ is
\begin{eqnarray}\label{boun1}
&&r=0~~\rightarrow~~x=0,\\\nonumber
&&\rho=\lambda~~\rightarrow~~y(0)=1.
\end{eqnarray}
Consequently, an additional condition must be introduced in order to
maintain the condition of Eq. (\ref{boun1}) simultaneously:
\begin{eqnarray*}
y'(0)=0~.
\end{eqnarray*}
In other words, the boundary conditions are as follows
\begin{equation*}
y(0)=1 ~~, ~~ y'(0)=0.
\end{equation*}
The values of $m$, which are physically interesting, lie in the interval
$[0, 5]$. The main difficulty in analysing this type of equation is
the singularity behaviour occurring at $x = 0$.\\
As it has been mentioned in the literature review, the solutions of the Lane-Emden equation could be exact only for $m=0, 1$ and 5. For the other values of $m$, the Lane-Emden equation is to be integrated numerically \citep{Horedt}.
Thus, we decided to present a new and
efficient technique to solve it numerically for $m=0, 1, 1.5, 2, 2.5, 3, 4$ and $5$.
\subsection{ \textbf{Methods have been proposed to solve Lane-Emden type equation}}
Many problems in science and engineering arise in unbounded domains. Different numerical methods have been proposed for solving problems on unbounded domains\citep{kumar2014analytical, rashidi2014free, kazem2011new, parand2011collocation, parand2011improved, parand2010modified, parand2013sinc, parand2010solution}. \\
Recently, many analytical methods have been used to solve Lane-Emden equations, the main difficulty arises in the singularity of the equation at $x=0$. Currently, most techniques in use for handling the Lane-Emden-type problems are based on either series solutions or perturbation techniques.\\
Bender et al. \citep{Bender}, proposed a new perturbation technique based on an artificial parameter $\delta$, the method is often called $\delta-$method.\\
Mendelzweig and Tabkin. \citep{Mendelzweig} used quasi-linearization approach to solve Lane-Emden equation. This method approximates the solution of a non-linear differential equation by treating the non-linear terms as a perturbation about the linear ones, and unlike perturbation theories is not based on the existence of some kind of a small parameter. He showed that the quasi-linearization method gives excellent results when applied to different non-linear ordinary differential equations in physics, such as the Blasius, Duffing, Lane-Emden and Thomas-Fermi equations.\\
Shawagfeh \citep{Shawagfeh} applied a non-perturbative approximate analytical solution for the Lane-Emden equation using the Adomian decomposition method, his solution was in the form of a power series. He used Pad\'e approximants method to accelerate the convergence of the power series.\\
In \citep{Wazwaz1}, Wazwaz employed the Adomian decomposition method with an alternate framework designed to overcome the difficulty of the singular point. It was applied to the differential equations of Lane-Emden type.\\
Later on \citep{Wazwaz2} he used the modified decomposition method for solving analytical treatment of non-linear differential equations such as Lane-Emden equations. The modified method accelerates the rapid convergence of the series solutions, dramatically reduces the size of work, and provides the solution by using few iterations only without any need to the so-called Adomian polynomials.\\
Liao \citep{LiaO} provided a reliable, easy-to-use analytical algorithm for Lane-Emden type equations. This algorithm logically contained the well-known Adomian decompositions method. Different from all other analytical techniques, this algorithm itself provides us with a convenient way to adjust convergence regions even without Pad\'e technique.\\
J.H. He \citep{He} employed Ritz's method to obtain an analytical solution of the problem. By the semi-inverse method, a variational principle is obtained for the Lane-Emden equation, which he gave much numerical convenience when applied to finite element methods or Ritz method.\\
Parand et al. \citep{Parand1,Parand2,Parand3} presented some numerical techniques to solve higher ordinary differential equations such as Lane-Emden. Their approach was based on a rational Chebyshev and rational Legendre tau method. They presented the derivative and product operational matrices of rational Chebyshev and rational Legendre functions.\\
These matrices together with the tau method were utilized to reduce the solution of these physical problems to the solutions of systems of algebraic equations. Also, in some paper's \citep{Parand4,Parand5,Parand6,Parand7,Parand8,Parand9},Parand et al. applied pseudo-spectral method based on rational Legendre functions, Sinc collocation method, Lagrangian method based on modified generalized Laguerre function, Hermite function collocation method ,Bessel function collocation method and meshless collocation method based on Radial basis function (RBFs) to solving the Lane-Emden type equations.\\
Ramos \citep{Ramos1,Ramos2,Ramos3,Ramos4} solved Lane-Emden equations through different methods. In \citep{Ramos1} he presented linearzation methods for singular initial value problems in second order ordinary differential equations such as Lane-Emden. These methods result in linear constant-coefficient ordinary differential equations which can be integrated analytically, thus they yield piecewise analytical solutions and globally smooth solutions \citep{Ramos2}. Later, he developed piecewise-adaptive decomposition methods for the solution of non-linear ordinary differential equations. Piecewise-decomposition methods provide series solutions in intervals which are subject to continuity conditions at the end points of each interval and their adoption is based on the use of either a fixed number of approximants and a variable step size, a variable number of approximants and a fixed step size or a variable number of approximants and a variable step size.\\
In \citep{Ramos3}, series solutions of the Lane-Emden equations based on either a volterra integral equations formulation or the expansion of the dependent variable in the original ordinary differential equations are presented and compared with series solutions obtained by means of integral or differential equations based on a transformation of the dependent variables.\\
Yousefi \citep{Yousefi} used integral operator and converted Lane-Emden equations to integral equations and then applied Legendre Wavelet approximations. In this work properties of Legendre wavelet together with the Gaussian integration method were utilized to reduce the integral equations to the solution of algebraic equations. By his method, the equations was formulated on [0, 1].\\
Chowdhury and Hashim \citep{Hashim} obtained analytical solutions of the generalized Emden-Fowler type equations in the second order ordinary differential equations by homotopy-perturbation method (HPM). This method is a coupling of the perturbation method and the homotopy method. The main feature of the HPM \citep{Dehghan2} is that it deforms a difficult problem into a set of problems which are easier to solve. HPM yields solution in convergent series forms with easily computable terms.\\
Aslanov \citep{Aslanov} constructed~a recurrence relation for the components of the approximate solution and investigated the convergence conditions for the Emden-Fowler type of equations. He improved the previous results on the convergence radius of the series solution.\\
Dehghan and Shakeri \citep{Dehghan3} investigated Lane-Emden equations using the variational iteration method 
and showed the efficiency and applicability of their procedure for solving this equations. Their technique does not require any discretization, linearization or small perturbations and therefore reduces the volume of computations.\\
Bataineh et al. \citep{Bataineh} obtained analytical solutions of singular initial value problems (IVPs) of the Emden-Fowler type by the homotopy analysis method (HAM). There solutions contained an auxiliary parameter which provided a convenient way of controlling the convergence region of the series solutions. It was shown that the solutions obtained by the Adomian decomposition method (ADM) and the Homotopy-perturbation method (HPM) are only special cases of the HAM solutions.\\
Marzban et al. \citep{Marzban} used a method based upon hybrid function approixmations. He used properties of hybrid of block-pulse functions and Lagrange interpolating polynomials together with the operational integration matrix for solving non-linear second-order, initial value problems and the Lane-Emden equations.\\
In \citep{Singh} Singh et al. used the modified Homotopy analysis method for solving the Lane-Emden-equations and White-Dwarf equation.\\
Adibi and Rismani in \citep{Adibi} proposed the approximate solutions of singular IVPs of the Lane-Emden type in second-order ordinary differential equations by improved Legendre-spectral method. The Legendre-Gauss point used as collocation nodes and Lagrange interpolation is employed in the Volterra term.\\
In \citep{Vanani} Karimi vanani and Aminataei provide a numerical method which produces an approximate polynomial solution for solving Lane-Emden equations as singular initial value problem. They are initially, used an integral operator and convert Lane-Emden equations into integral equations. Then, convert the acquired integral equations into a power series and finally, transforming the power series into pad\'e series form.\\
Kaur et al. \citep{Kaur}, obtained the Haar wavelet approximate solutions for the generalized Lane-Emden equations. These method was based on the quasi-linearization approximation and replacement of an unknown function through a truncated series of Haar wavelet series of the function.\\
so, the other researchers trying to solving the Lane-Emden type equtions with several methods, For example, Y{\i}ld{\i}r{\i}m and \"{O}zi\c{s} \citep{Ozis1,Ozis2} by using HPM and VIM methods, Benko et al. \citep{Benko} by using Nystr\"{o}m method, Iqbal and javad \citep{IQbal} by using Optimal HAM, Boubaker and Van Gorder \citep{Boubaker} by using boubaker polynomials expansion scheme, Da\c{s}c{\i}o\u{g}lu and Yaslan \citep{Dascioglu} by using Chebyshev collocation method, Y\"{u}zba\c{s}{\i} \citep{Yuzbasi1,Yuzbasi2} by using Bessel matrix and improved Bessel collocation method, Boyd \citep{Boyd} by using Chebyshev spectral method, Bharwy and Alofi \citep{Alofi} by using Jacobi-Gauss collocation method, Pandey et al. \citep{Pandey1,Pandey2} by using Legendre and Brenstein operation matrix, Rismani and monfared \citep{Rismani} by using Modified Legendre spectral method, Nazari-Golshan et al. \citep{Golshan} by using Homotopy perturbation with Fourier transform, Doha et al. \citep{Doha} by using second kind Chebyshev operation matrix algorithm, Carunto and bota \citep{bota} by using Squared reminder minimization method, Mall and Chakaraverty \citep{Mall} by using Chebyshev Neural Network based model, G\"{u}rb\"{u}z and sezer \citep{Gurbuz} by using Laguerre polynomial and Kazemi-Nasab et al. \citep{Kazemi} by using Chebyshev wavelet finite difference method.\\
In this paper, we attempt to introduce a new method, based on $RBF-DQ$ for solving non-linear ODEs. \\
The organization of the paper is as follows:\\
In section 2, we briefly explain on the RBF, DQ and RBF-DQ methods. In section 3 we applying the RBF-DQ method and Solving some Lane-Emden type equations.Then, a comparison is made with the existing methods in the literature. section 4 is devoted to conclusions.
\section{\textbf{Development of radial basis function- based differential quadrature (RBF- DQ)method }}
In this section, we will show in detail the RBF-DQ method. The RBF-differential quadrature (DQ) method is an interpolation technique in which the radial basis functions are used as basis functions.So, it can be easily applied to the linear and non-linear problems. In the following, we will show the details of RBF-DQ method step by step. 
\subsection{\textbf{Radial basis functions (RBFs)}}
The interpolation of a given set of points is an important problem especially in higher dimensional domains. Although polynomials (e.g. Chebyshev and Legendre) are very powerful tools for interpolating of a set of points in one-dimensional domains, the use of these functions are not efficient in higher dimensional or irregular domains. While applying these functions, the points in the domain of the problem should be chosen in a special form, which is very limiting when the interpolation of a scattered set of points is needed. Radial basis functions (RBFs) are very efficient instruments for interpolating a scattered set of points, which have been used in the last 20 years. The use of radial basis functions collocation method for solving partial differential equations has some advantages over mesh dependent methods, such as finite-difference methods, finite element methods, spectral elements, finite volume methods and boundary element methods. Since a large portion of the computational time is spent for providing a suitable mesh on the domain of the problem in mesh-dependent methods, the meshfree methods have an auxiliary role in the numerical solution of partial differential equations. In recent years, to avoid the mesh generation, meshfree methods have attracted the attention of researchers.\par
In a meshless method, a set of scattered nodes is used instead of meshing the domain of the problem. RBFs were first studied by Roland Hardy, an Iowa State geodesist, in 1968. This method allows scattered data to be easily used in computations. An extensive study of interpolation methods available at the time was conducted by Franke \citep{Franke}, who concluded that RBF interpolations were evaluated as the most accurate techniques. The theory of RBFs originated as a means to prepare a smooth interpolation of a discrete set of data points. The concept of using RBFs for solving differential equations (DEs) was first introduced
by Kansa \citep{Kansa1,Kansa2}, who directly collocated the radial basis functions for the approximate solution of DEs. However, in recent years RBFs have been extensively researched and applied in a wider range of analysis. Partial differential equations (PDEs) and ordinary differential equations (ODEs) have been solved using RBFs in recent work \citep{Power, Mai-Duy, ParandAbbasbandy1, Islam, kazem1, kazem2, kazem3, kazem4, kazem5}.\par
A well-known space RBF $\varphi(\parallel X$-$X_i\parallel) : \mathbb{R}^{+}\longrightarrow \mathbb{R}$ depends on the separation of a field point $ X \in \mathbb{R}^d$ and the data centers$ X_i$ , for$ i = 1, 2, \cdots, N,$ and N data points. The interpolants are classed as radial function due to their spherical symmetry around centres $X_i$, where $\parallel . \parallel$ is the well-known Euclidean norm. The most known space RBFs are listed in Table \ref{TabRBF}, where r = $\varphi(\parallel X$-$X_i\parallel)$ and $c$ is a free positive parameter, often referred to as the shape parameter, to be specified by the user. The shape parameter $c$ within the Gaussian and multiquadric RBFs requires fine tuning and can dramatically alter the quality of the interpolation. Too large or too small shape parameter makes the GA too flat or too peaked.
Despite many research works, which are done to find algorithms for selecting the optimum values of $\epsilon$ \citep{Rippa, Cheng, Carlson, Tarwater, Fasshauer} the optimal choice of shape parameter is an open problem, which is still under intensive investigation. One of the most powerful interpolation methods with analytic d-dimensional test function is the space RBFs method, based on Gaussian (GA) basis function
\begin{equation*}
\varphi(r)=e^{-c^{2}r^2},~~~~~c>0.
\end{equation*}
In the cases of inverse quadratic, inverse multiquadric (IMQ) and Gaussian (GA), the coefficient matrix of RBFs interpolating is positive definite and, for multiquadric (MQ), it has one positive eigenvalue and the remaining ones are all negative \citep{Powel}.
The interested reader is referred to the recent books and the paper written by Buhmann \citep{Buhmann1, Buhmann2} and Wendland \citep{Wendland} for more basic details about RBFs, compactly and globally supported and convergence rate of the radial basis functions.
\begin{table}[htbp]
\centering 
\begin{tabular}{l l l} 
\hline\hline
$Name~of~functions$ &$Definition$ \\
\hline
$Multiquadrics (MQ)$ &$\sqrt{ r^{2}+c^{2}}$\\
$Inverse ~multiquadrics (IMQ)$ &$\frac{1}{\sqrt{ r^{2}+c^{2}}}$\\
$Gaussian (GA)$ &$\exp^{-c^{2}r^2}$\\
$Inverse ~quadrics$&$\frac{1}{r^2+c^2}$\\
$Thin~ Plate ~Splines (TPS)$&$(−1)^{k+1}+r^{2k} log(r) $\\
$Hyperbolic ~secant (sech)$&$sech(c\sqrt{r})$\\
\hline
\end{tabular}
\caption{Some well-known RBFs$(r=\parallel X$-$X_i\parallel^2=r_i)$, $c>0$}.
\label{TabRBF}
\end{table}
\subsection{\textbf{The differential quadrature (DQ) method}}
Solving differential equations (ordinary and partial) is one of the most important ways engineers, physicists and applied mathematicians use to tackle a practical problem. 
The differential quadrature (DQ) method was presented by R.E.Bellman and his associations in the early 1970's.\par
The DQ method is a numerical discretization technique to approximate of derivatives.This method was initiated from the idea of conventional integral quadrature.
\begin{equation}\label{eqINT}
\int_{a}^{b} f(x)~dx =w_1f_1+w_2f_2+\cdots+w_nf_n=\sum_{k=1}^{N} w_kf_k
\end{equation}
Where $ w_1,w_2,\cdots w_n $ are weighting coefficients, $f_1,f_2,.\cdots ,f_n$ are the functional values at the discrete points $a=x_1,x_2,\cdots,x_n=b$.\par
Following the idea of integral quadrature Eq. (\ref{eqINT}), Bellman at all (1972) \citep{bellman} suggested that the first order derivative of the function $f^{(1)}(x)$ with respect to x at points $x_i$,is approximated by a linear combination of all functional values in the whole domain by
\begin{equation}\label{eqDQ1}
f^{(1)}_x(x_i)=\frac{\partial f}{\partial x}|_{x_i}=\sum_{j=1}^{N} w^{(1)}_{ij}f(x_j)~~~~~ i =1,2,\cdots,N
\end{equation}
Where $f(x_j) $ represents the functional value at a grid point $x_j $ and $w^{(1)}_{ij}$ represented the weighting coefficients, and $N$ is the number of points the whole domain .\par
So, the m-th order derivative with respect to x at a point $ x_i$ can be approximated by DQ as 
\begin{equation}\label{eqDQm}
f^{(m)}_x(x_i)=\frac{\partial^{m} f}{\partial x^m}|_{x_i}=\sum_{j=1}^{N} w^{(m)}_{ij}f(x_j)~~~~~ i =1,2,\cdots,N
\end{equation}
where $x_j$ are the discrete points in the domain, $f(x_j)$ and $w^{(m)}_{ij}$ are the function values at these points and the
related weighting coefficients.Obviously, the key procedure in the DQ method is the determination of the weighting coefficients $w^{(m)}_{ij}$. It has been shown by Shu \citep{shuDQ, zongDQ} that the weighting coefficients can be easily computed under the analysis of a linear vector space and the analysis of a function approximation. 
\subsection{\textbf{The RBF-DQ method}}
According to DQ definition, after domain discretization, the n-th order derivatives of a function f(x) with respect to x, $f^{(n)}_x$, at the i-th node $x_i$ can be written as 
\begin{equation}\label{eq16}
f^{(n)}_x(x_i)=\sum_{j=1}^{N} w^{(n)}_{ij}f(x_j)
\end{equation}
In the usual procedure of RBF-DQ method, we combine this method with the radial basis function (RBF) based interpolation of the function f. If $\varphi$ is a RBF and we consider the center point $x_i$, we may write the interpolant to the function f as:\\
\begin{equation}\label{eq17}
f(x)=\sum_{k=1}^{N} \lambda_k \varphi_k(x)
\end{equation}
where $\lambda_k$ is the coefficient for $\varphi_k(x)$ and $\varphi_k(x)$ is a RBF.\par
Since the RBF values and their n-th order derivatives are known at each node, the unknown weighting coefficients can be
determined by solving a linear system of equations:
\begin{equation}\label{eq18}
\frac{\partial^{(n)} \varphi_k(x_i)}{\partial x^{(n)}} =\sum_{j=1}^{N} w^{(n)}_{ij}\varphi_k(x_j).\\
\end{equation}
The above equation can be further put in the matrix form 
$$\underbrace{
\begin{pmatrix}
\frac{\partial^{(n)} \varphi_1(x_i)}{\partial x^{(n)}} &\cr 
\frac{\partial^{(n)} \varphi_2(x_i)}{\partial x^{(n)}} &\cr \vdots \cr
\frac{\partial^{(n)} \varphi_k(x_i)}{\partial x^{(n)}}\cr
\end{pmatrix}
}
_{\frac{\partial \overrightarrow{\varphi}(x_i)}{\partial x}}
=
\underbrace{
\begin{pmatrix}
\varphi_1(x_1)&\varphi_1(x_2)&\cdots &\varphi_1(x_n)\cr
\varphi_2(x_1)&\varphi_2(x_2)&\cdots &\varphi_2(x_n)\cr \vdots &\vdots &\ddots &\vdots \cr
\varphi_n(x_1)&\varphi_n(x_2)&\cdots &\varphi_n(x_n)\cr
\end{pmatrix}
}_{(A)}
\underbrace{
\begin{pmatrix}
w^{(n)}_{i1}&\cr
w^{(n)}_{i2}&\cr \vdots \cr
w^{(n)}_{in}\cr
\end{pmatrix}
}_{\overrightarrow{w}_{i}}
.$$
According to the theory of RBF approximation, we can write,
\begin{equation}\label{eq20}
\overrightarrow{w}^{n}_{i}=
{\begin{pmatrix}
A
\end{pmatrix}}^{-1} \frac{\partial^{(n)} \overrightarrow{\varphi}(x_i)}{\partial x^{(n)}}
\end{equation}
Note that$ \varphi_k(x_j)=\varphi(\parallel x_j-x_k \parallel)$. With solving the system Eq.\eqref{eq20}, we obtain the
unknown weighting coefficients vector $w_i$.
\section{\textbf{Applying the RBF-DQ method for solving the Lane-Emden equations} }
In this section, we explain function approximation based on RBF-DQ method to solve of Lane-Emden type equations. In general the Lane-Emden type equation are formulated as Eq.\eqref{Lane-Emden} with initial conditions Eq.\eqref{conditions} \\
According to Eq. (\ref{eq16}), an approximation derivative of the first and second order of function $f(x)$ with respect to x, $f^{(1)}_{x}, f^{(2)}_{x}$, at the ith nod $x_i$, can be writen\\
\begin{equation}\label{eq31}
f^{(1)}_x(x_i)=\sum_{j=1}^{N} w^{(1)}_{ij}f(x_j)=\sum_{j=1}^{N} w^{(1)}_{ij}f_j,
\end{equation}
\begin{equation}\label{eq32}
f^{(2)}_x(x_i)=\sum_{j=1}^{N} w^{(2)}_{ij}f(x_j)=\sum_{j=1}^{N} w^{(1)}_{ij}f_j.
\end{equation}
Where $w^{(1)}_{ij} ~and ~w^{(2)}_{ij}$ are the weighted coefficient for $f(x_j)$ and $x_i$ obtained N-nodes as following equation\\
\begin{equation}\label{eq33}
x_i=\frac{L}{2}(1-cos(\frac{i-1}{N-1})\pi))~~~~i=1,2,\cdots,N.
\end{equation}
Where L is the upper bound of domain problem.\par
Now, with respect the initial value of problem, we have
\begin{equation}\label{eq34}
\begin{cases}
f(x_0)=f_1=1\\
\sum_{i=1}^{N} w^{(1)}_{1,i}f_i=0.
\end{cases}
\end{equation}
Then, to apply the collocation mehod, we produced the residual function as following,
\begin{equation}\label{eq35}
Res(x)=x.\sum_{j=1}^{N} w^{(2)}_{i,j}f_j+2\sum_{j=1}^{N} w^{(1)}_{i,j}f_j+x(f_i)^m~~~~~i=1,2,\cdots,N.
\end{equation}
By using $x_i$ and residual function, we obtained $N$ equation, so by solving these equations we obtained the $f_i$s
\begin{equation}\label{eq36}
Res(x_i)=0~~~~~~~i=1,2,\cdots,N.
\end{equation}
\subsection{\textbf{Solving some Lane-Emden type equations}}
\subsubsection{\textbf{Example 1 (The standard Lane-Emden equation)}}
\begin{figure}[htbp]
\center
\includegraphics[scale=.45]{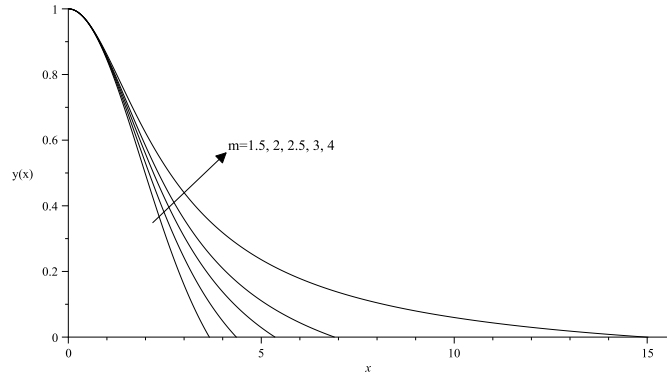}\\
\caption{The obtained graphs of solutions of Lane-Emden standard equations for m=1.5, 2, 2.5 ,3 , 4 }\label{Lane-Emden-standard}
\end{figure}

For $f(x,y)=y^{m}$, $A=1$ and $B=0$, Eq.\eqref{Lane-Emden} is the standard Lane-Emden equation, which was used to model the thermal behaviour of a spherical cloud of gas acting under the mutual attraction of its molecules and subject to the classical laws of thermodynamic \citep{Shawagfeh}
\begin{equation}\label{Ex1}
y''(x) + \frac {2}{x} y'(x) + y^{m}(x) = 0, ~~ x>0
\end{equation}
with conditions:
\begin{eqnarray*}
y(0)=1,\\
y'(0)=0,
\end{eqnarray*}
where $m\geq0$ is constant. For $m=0,1$ and $5$ Eq.\eqref{Ex1} has the exact solution, respectively:
\begin{equation*}
y(x)=1-\frac{1}{3!}x^{2},~~~~ y(x)=\frac{\sin(x) }{x}, ~~~~y(x)=\Bigg(1+\frac{x^2}{3}\Bigg)^{-\frac{1}{2}}~.
\end{equation*}
In other cases, there is not any exact analytical solution. Therefore, we apply the RBF-DQ method to solve the standard Lane-Emden
Eq. (\ref{Ex1}), for $m=1.5$, $2$, $2.5$, $3$ and $4$. In this example, we construct the residual function as follows:
\begin{equation*}
Res(x)=x.\sum_{j=1}^{N} w^{(2)}_{i,j}u_j+2\sum_{j=1}^{N} w^{(1)}_{i,j}u_j+x(f_i)^m~~~~~i=1,2,\cdots,N.
\end{equation*}
Therefore, to obtain the coefficient $w^{(1)}_{i,j}~and~w^{(2)}_{i,j}$, $Res(x)$ is equalized to zero at $N$ collocation point. By solving this set of non-linear algebraic equations, we can find the approximating function $f_i(x)$.\\
Table \ref{table:m} shows the comparison of the first zero of standard Lane-Emden equations, from the present method and exact values m given by Hordet \citep{Horedt} and for m=1.5 , 2, 2.5, 3 and 4, respectively and Fig.\ref{Lane-Emden-standard} represents the graphs of the standard Lane-Emden equations for m=1.5, 2, 2.5, 3 and 4. \\
Tables $3-6$ show the obtained $y(x)$ for the standard Lane-Emden equations for m=1.5, 2, 2.5, 3 and 4, respectively, by the proposed RBF-DQ method in this paper, and some well-known method in other articles. Also these tables show the residual function $Res(x)$ in some points of the interval $[0,\infty)$.
These tables demonstrate that the present method has a good accuracy.
\begin{table}[htbp]

\centering 
\begin{tabular}{l l| l l l} 
\hline\hline 
$m$ &N & present method & Horedt \citep{Horedt} \\ [0.5ex] 
\hline 
$1.5$ & $30$& $3.6537537342$ & $3.65375374$ \\ 
$2$ & $40$& $4.3528745959$ & $4.35287460$ \\
$2.5$ & $30$& $5.3552754590$ & $5.35527546$ \\
$3$ & $30$& $6.8968486191$ & $6.89684862$ \\
$4$ & $60$& $14.971538050$ & $14.9715463$ \\[1ex] 
\hline 
\end{tabular}
\caption{Obtained first zeros of standard Lane-Emden equations, by the present method for several $m$} 
\label{table:m} 
\end{table}
\begin{table}[htbp]
\centering 
\begin{tabular}{l l l l l} 
\hline\hline 
$x$ &Horedt \citep{Horedt} & present Method & $Res(x)$ \\ [0.5ex] 
\hline 
$0.00$ & $1.0000000$ & $1.0000000000$ & $0.00000$ \\ 
$0.10$ & $0.9983346$ & $0.9983345826$ & $8.891e-10$ \\
$0.50$ & $0.9591039$ & $0.9591038569$ & $1.322e-9$\\
$1.00$ & $0.8451698$ & $0.8451697549$ & $1.016e-7$\\
$3.00$ & $0.1588576$ & $0.1588576082$ & $4.452e-7$\\
$3.60$ & $0.1109099e-1$ & $0.1109099415e-1$ & $6.433e-6$\\
$3.65$ & $0.7639242e-3$ & $0.7639240088e-3$ & $4.988e-6$\\ [1ex] 
\hline 
\end{tabular}
\caption{Obtained values of $y(x)$ for Lane-Emden standard $m=1.5$ by present method ($c=1$ and $N$=30).} 
\label{table:m=1.5} 
\end{table}
\begin{table}[htbp]
\centering 
\begin{tabular}{l l l l l} 
\hline\hline 
$x$ & Horedt \citep{Horedt} & present Method & $Res(x)$\\ [0.5ex] 
\hline 
$0.0$ & $1.0000000$ & $1.00000000000$ & $0.0000000$ \\ 
$0.1$ & $0.9983350$ & $0.99833499854$ & $6.9849e-12$ \\
$0.5$ & $0.9983350$ & $0.95935271580$ & $8.9006e-11$ \\
$3.0$ & $2.418241e-1$ & $0.24182408305$ & $3.9602e-11$ \\
$4.3$ & $6.810943e-3$ & $6.81094327419e-3$ & $1.1602e-10$ \\
$4.35$ & $3.660302e-4$ & $3.66030179364e-4$ & $1.5327e-10$\\[1ex] 
\hline 
\end{tabular}
\caption{Obtained values of $y(x)$ for standard Lane-Emden $m=2$ by the present method (with $c=1$ and $N=30$).} 
\label{table:m=2} 
\end{table}
\begin{table}[htbp]
\centering 
\begin{tabular}{l l l l l} 
\hline\hline 
$x$ & Horedt \citep{Horedt} & present method & $Res(x)$\\ [0.5ex] 
\hline 
$0.0$ & $1.000000$ & $1.00000000000$ & $0.0000$ \\ 
$0.1$ & $9.983354e-1$ &$9.98335414193e-1$ & $8.8589e-10$ \\
$0.5$ & $9.595978e-1$ & $9.59597754452e-1$ & $6.7729e-9$ \\
$1.0$ & $8.519442e-1$ & $8.51944199404e-1$ & $4.0086e-8$ \\
$4.0$ & $1.376807e-1$ & $1.37680732920e-1$ & $1.2091e-7$ \\
$5.0$ & $2.901919e-2$ & $2.90191866504e-2$ & $7.4321e-9$ \\
$5.3$ & $4.259544e-3$ & $4.25954353341e-3$ & $3.4189e-9$ \\
$5.355$ & $2.100894e-5$ & $2.10089389496e-5$ & $1.2909e-7$\\[1ex] 
\hline 
\end{tabular}
\caption{Obtained values of $y(x)$ for standard Lane-Emden $m=2.5$ by the present method (with $c=1$ and $N=30$ ).} 
\label{table:m=2.5} 
\end{table}
\begin{table}[htbp]
\centering 
\begin{tabular}{l l l l l} 
\hline\hline 
$x$ & Horedt \citep{Horedt} & present method & $Res(x)$ \\ [0.5ex] 
\hline 
$0.0$ & $1.00000000$ & $1.00000000000$ & $0.0000$ \\ 
$0.1$ & $9.983358e-1$ & $9.9833583002e-1$ & $3.7571e-8$ \\
$0.5$ & $9.598391e-1$ & $9.5983906655e-1$ & $6.8415e-8$ \\
$1.0$ & $8.550576e-1$ & $8.5505753413e-1$ & $3.3963e-6$ \\
$5.0$ & $1.108198e-1$ & $1.1081993474e-2$ & $2.8895e-5$ \\
$6.0$ & $4.373798e-2$ & $4.3737983486e-2$ & $1.1100e-6$ \\
$6.8$ & $4.167789e-3$ & $4.1677893542e-3$ & $1.1411e-8$ \\
$6.896$ & $3.601115e-5$ & $3.6011134911e-5$ & $1.5511e-7$ \\[1ex] 
\hline 
\end{tabular}
\caption{Obtained values of $y(x)$ for standard Lane-Emden $m=3$ by present the method (with $c=1$ and $N=30$).} 
\label{table:m=3} 
\end{table}
\begin{table}[htbp]
\centering 
\begin{tabular}{l l l l l} 
\hline\hline 
$x$ & Horedt \citep{Horedt} & present method & $Res(x)$ \\ [0.5ex] 
\hline 
$0.0$ & $1.000000$ & $1.000000000$ & $0.0000$ \\ 
$0.1$ & $9.983367e-1$ & $9.983366e-1$ & $8.4117e-8$ \\
$0.2$ & $9.933862e-1$ & $9.933862e-1$ & $7.6620e-7$ \\
$0.5$ & $9.603109e-1$ & $9.603109e-1$ & $1.23006e-5$ \\
$1.0$ & $8.608138e-1$ & $8.608144e-1$ & $1.3033e-4$ \\
$5.0$ & $2.359227e-1$ & $2.352433e-1$ & $9.3449e-2$ \\
$10.0$ & $5.967274e-2$ & $5.965197e-2$ & $4.4052e-3$ \\
$14.0$ & $8.330527e-3$ & $8.330447e-3$ & $2.5998e-5$ \\
$14.9$ & $5.764189e-4$ & $5.763524e-4$ & $1.1584e-6$ \\[1ex] 
\hline 
\end{tabular}
\caption{Obtained values of $y(x)$ for standard Lane-Emden $m=4$ by the present method (with $c=1$ and $N=30$)} 
\label{table:m=4} 
\end{table}\\\\\\\\\\\\\\\\\\\\\\\\\\
\subsubsection{\textbf{Example 2 (The isothermal gas spheres equation)}}
For $f(x,y)=e^{y}$, $A=0$ and $B=0$, (\ref{Lane-Emden}) is the isothermal gas sphere equation
\begin{equation}\label{Ex2}
y''(x) + \frac {2}{x} y'(x) + e^{y(x)} = 0, ~~ x>0,
\end{equation}
with conditions
\begin{eqnarray*}
y(0)=0,\\
y'(0)=0.
\end{eqnarray*}
This model can be used to view the isothermal gas sphere. For reading the thorough discussion of Eq. (\ref{Ex2}), see \citep{Davis}. This equation has been solved by some researchers, for example \citep{Wazwaz1, Parand7, MSH} by ADM, Hermit collocation method and HPM, respectively.\\
\par A series solution obtained by Wazwaz \citep{Wazwaz1} by using ADM and series expansion:
\begin{eqnarray}\label{Eq2}
y(x)\simeq - \frac{1}{6}x^{2}+\frac{1}{5. 4!}x^4 - \frac{8}{21. 6!}x^{6}+\frac{122}{81. 8!}x^{8}-\frac{61.67}{495.10!}x^{10}.
\end{eqnarray}
The residual function for Eq. (\ref{Ex2}) is as follows:
\begin{table}[htbp]
\scriptsize\centering 
\begin{tabular}{l l l l l l} 
\hline\hline 
$x$ &Peresent Method & BFC \citep{Parand8} & ADM\citep{Wazwaz1} & HFC\citep{Parand7} & Error\\ [0.5ex] 
\hline 
$0.0$ & $~~0.0000000$ & $~~0.0000000$ & $~~0.0000000$ & $~~0.0000000$ & $0.0000000$\\ 
$0.1$ & $-0.00166583$ & $-0.00166583$ & $-0.0016658$ & $-0.0016664$ & $1.13e-14$ \\ 
$0.2$ & $-0.00665336$ & $-0.00665336$ & $-0.0066533$ & $-0.0066539$ &$1.21e-14$ \\ 
$0.5$ & $-0.04115395$ & $-0.04115395$ & $-0.0411539$ & $-0.0411545$ & $5.25e-14$ \\
$1.0$ & $-0.15882767$ & $-0.15882767$ & $-0.1588273$ & $-0.1588281$ & $2.46e-13$ \\ 
$1.5$ & $-0.33801942$ & $-0.33801942$ & $-0.3380131$ & $-0.3380198$ & $5.16e-13$ \\ 
$2.0$ & $-0.55982300$ & $-0.55982300$ & $-0.5599626$ & $-0.5598233$ & $5.48e-13$ \\ 
$2.5$ & $-0.80634087$ & $-0.80634087$ & $-0.8100196$ & $-0.8063410$ & $3.26e-13$ \\[1ex] 
\hline 
\end{tabular}
\caption{The comparison of the values $y(x)$ of Lane-Emden equation obtained by the present method with $c=1$ and $N=30$ and \citep{Wazwaz1, Parand7, Parand8}and Bessel in Example 2.} 
\label{TabEX2} 
\end{table}\\
\begin{figure}[htbp]
\centering
\includegraphics[width=2in]{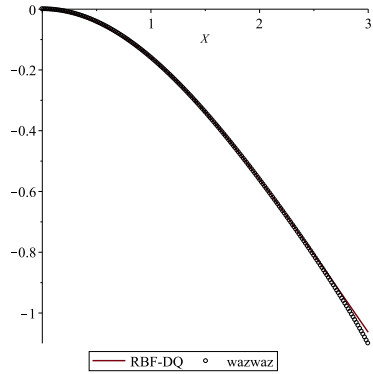}\\
\caption{Graph of isothermal gas sphere equation in comparison with wazwaz solution \citep{Wazwaz1}}
\label{figex2}
\end{figure}
\begin{equation*}
Res(x)=x.\sum_{j=1}^{N} w^{(2)}_{i,j}f_j+2\sum_{j=1}^{N} w^{(1)}_{i,j}f_j+xe^{(f_i)}~~~~~i=1,2,\cdots,N.
\end{equation*}
By solving the set of equations obtained from $N$ equations $Res(x)=0$ at $N$ collocation points, we have the approximating function $f_i$.\\
Table \ref{TabEX2} shows the comparison of $y(x)$ obtained by the proposed method in this paper and \citep{Wazwaz1, Parand7, Parand8}.\\
The resulting graph of the isothermal gas spheres equation in comparison to the presents method and those obtained by wazwaz \citep{Wazwaz1} is shown in Fig. \ref{figex2}.
\subsubsection{\textbf{Example 3}}
For $f(x,y)=\sinh(y)$, $A=1$ and $B=0$, Eq.\eqref{Lane-Emden} becomes
\begin{equation}\label{Ex3}
y''(x) + \frac {2}{x} y'(x) + \sinh(y(x)) = 0, \texttt{ } x>0,
\end{equation}
with conditions
\begin{eqnarray*}
y(0)=1,\\
y'(0)=0.
\end{eqnarray*}
This equations has been solved by some researchers, for example \citep{Wazwaz1, Parand7, Parand8}by ADM , Hermit and Bessel orthogonal functions collocation method respectively. A series solution obtained in \citep{Wazwaz1}, by using ADM, is:
\begin{eqnarray}\label{Eq3}
&&y(x)\simeq 1-\frac{(e^{2}-1)x^{2}}{12 e}+\frac{1}{480} \frac{(e^{4}-1)x^{4}}{e^{2}}-\frac{1}{30240} \frac{(2 e^{6}+3e^{2}-3e^{4}-2)x^{6}}{e^{3}}\\\nonumber
&&~~~~~~~~+\frac{1}{26127360} \frac{(61e^{8}-10e^{6}+10 e^{2}-61)x^{8}}{e^{4}}.
\end{eqnarray}
We apply RBF-DQ method to solve Eq.\eqref{Ex3}, therefore we construct the residual function as follows:
\begin{equation*}
Res(x)=x.\sum_{j=1}^{N} w^{(2)}_{i,j}f_j+2\sum_{j=1}^{N} w^{(1)}_{i,j}f_j+x\sinh(f_i)~~~~~i=1,2,\cdots,N. 
\end{equation*}
To obtain coefficients $w^{(1)}_{i,j}~and~w^{(2)}_{i,j}$, we put $Res(x)=0$ at the $N$ collocation points, and, then, we construct a set of $N$ equations. By solving this set, we have the approximating function $f_i$.
\begin{table}[htbp]
\scriptsize\centering 
\begin{tabular}{l l l l l l} 
\hline\hline 
$x$ & present Method & BFC\citep{Parand8} & ADM\citep{Wazwaz1} & HFC\citep{Parand7} & Res(x) \\ [0.5ex] 
\hline 
$0.0$ & $1.0000000$ &$1.0000000$ & $1.0000000$ & $1.0000000$ & $0.0000000$\\ 
$0.1$ & $0.998042841$ &$ 0.998042841$ & $0.9980428$ & $0.9981138$ & $3.49e-16$ \\
$0.2$ & $0.992189434$ &$ 0.992189434$ & $0.9921894$ & $0.9922758$ & $7.12e-16$ \\
$0.5$ & $0.951961092$ & $0.951961092$ & $0.9519611$ & $0.9520376$ & $2.11e-15$ \\
$1.0$ & $0.818242928$ & $0.818242928$ & $0.8182516$ & $0.8183047$ & $6.52e-14$ \\
$1.5$ & $0.625438763$ & $0.625438763$ & $0.6258916$ & $0.6254886$ & $1.69e-13$ \\
$2.0$ & $0.406622887$ & $0.406622887$ & $0.4136691$ & $0.4066479$ & $8.16e-14$ \\[1ex] 
\hline 
\end{tabular}
\caption{The comparison of the values of $y(x)$ of Lane-Emden equation obtained from the present method with $c=1$ and $N=30$ and \citep{Wazwaz1, Parand7, Parand8} in Example 3.} 
\label{TabEX3}
\end{table} 
\begin{figure}[htbp]
\centering
\includegraphics[width=2in]{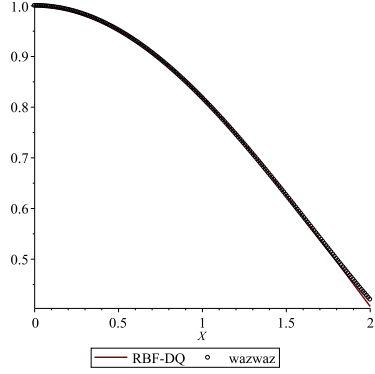}\\
\caption{Graph of equation of example 3 in comparing the presented method and Wazwaz solution \citep{Wazwaz1}}\label{figex3}
\end{figure}
\par
Table \ref{TabEX3} compares the $f(x)$ obtained by the proposed method in this paper and \citep{Wazwaz1, Parand7, Parand8}.\\
The resulting graph of the isothermal gas spheres equation in comparison to the presents method and those obtained by wazwaz \citep{Wazwaz1} is shown in Fig. \ref{figex3}.\\\\\\\\\\\\\
\subsubsection{\textbf{Example 4}}
For $f(x,y)=\sin(y(x))$, $A=1$ and $B=0$, Eq.\eqref{Lane-Emden} has the following form:
\begin{equation}\label{Ex4}
y''(x) + \frac {2}{x} y'(x) + \sin(y(x)) = 0,~~ x>0,
\end{equation}
with conditions
\begin{eqnarray*}
y(0)=1,\\
y'(0)=0.
\end{eqnarray*} 
A series solution obtained by Wazwaz \citep{Wazwaz1} by using ADM is:
\begin{eqnarray}\label{Eq4}
&&y(x)\simeq 1-\frac{1}{6}k x^{2}+\frac{1}{120}k l x^{4}+k(\frac{1}{3024}k^{2}-\frac{1}{5040}l^{2})x^{6}\\ \nonumber
&&~~~~~~~~+kl(-\frac{113}{3265920}k^{2}+\frac{1}{362880}l^{2})x^{8}\\\nonumber
&&~~~~~~~~+k(\frac{1781}{898128000}k^{2}l^{2}-\frac{1}{399168000}l^{4}-\frac{19}{2395080}k^{4})x^{10},
\end{eqnarray}
where $k=\sin(1)$ and $l=\cos(1)$.
Now we apply RBF-DQ method for Eq. (\ref{Eq4}) and construct the residual function as follows:
\begin{equation*}
Res(x)=x.\sum_{j=1}^{N} w^{(2)}_{i,j}f_j+2\sum_{j=1}^{N} w^{(1)}_{i,j}f_j+x\sin(f_i)~~~~~i=1,2,\cdots,N.
\end{equation*}
To obtain coefficients $w^{(1)}_{i,j}~and~w^{(2)}_{i,j}$, we put $Res(x)=0$ at the $N$ collocation points, and, then, we construct a set of $N$ equations. By solving this set, we have the approximating function $f_i$.
\begin{table}[htbp]
\scriptsize\centering 
\begin{tabular}{l l l l l l} 
\hline\hline 
$x$ &present Method & BFC\citep{Parand7} & ADM\citep{Wazwaz1} & HFC\citep{Parand7} & Res(x) \\ [0.5ex] 
$0.0$ & $1.0000000$ & $1.0000000$ & $1.0000000$ & $1.0000000$ & $0.0000000$\\ 
$0.1$ & $0.998597927$ & $0.998597927$ & $0.9985979$ & $0.9986051$ & $7.40e-16$ \\
$0.2$ & $0.994396264$ & $0.994396264$ & $0.9943962$ & $0.9944062$ & $1.49e-15$ \\
$0.5$ & $0.965177780$ & $0.965177780$ & $0.9651777$ & $0.9651881$ & $4.36e-15$ \\
$1.0$ & $0.863681125$ & $0.863681125$ & $0.8636811$ & $0.8636881$ & $1.30e-13$ \\
$1.5$ & $0.705045233$ & $0.705045233$ & $0.7050419$ & $0.7050524$ & $3.29e-13$ \\
$2.0$ & $0.506463627$ & $0.506463627$ & $0.5063720$ & $0.5064687$ & $1.51e-13$ \\[1ex] 
\hline 
\end{tabular}
\caption{The comparison of the values of $y(x)$ of Lane-Emden equation obtained from present method with $c=1$ and $N=30$ and \citep{Wazwaz1, Parand7, Parand8} in Example 4.} 
\label{TabEX4}
\end{table} 
\begin{figure}[htbp]
\centering
\includegraphics[width=2in]{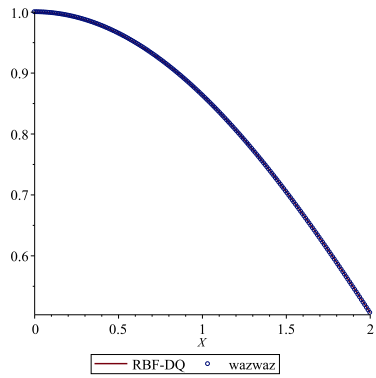}\\
\caption{Graph of equation of example 4 in comparing the presented method and Wazwaz solution \citep{Wazwaz1} }\label{figex4}
\end{figure}
\par
Table \ref{TabEX4} compares the $f(x)$ obtained by the proposed method in this paper and \citep{Parand7, Parand8,Wazwaz1} .\\
The resulting graph of the isothermal gas spheres equation in comparison to the presents method and those obtained by wazwaz \citep{Wazwaz1} is shown in fig.\ref{figex4}.
\subsubsection{\textbf{Example 5}}
For $f(x,y)=4 (2 e^{y(x)}+e^{\frac{y(x)}{2}})$, $A=0$ and $B=0$, Eq.\eqref{Lane-Emden} has the following form:
\begin{equation}\label{Eq5}
y''(x) + \frac {2}{x} y'(x) + 4 (2 e^{y(x)}+e^{\frac{y(x)}{2}}) = 0, ~~ x>0,
\end{equation}
with conditions:
\begin{eqnarray*}
y(0)=0,\\
y'(0)=0,
\end{eqnarray*}
which gives the following analytical solution:
\begin{equation} \label{Ex5}
y(x)= -2 \ln(1+x^{2}).
\end{equation}
This type of equation has been solved by \citep{Parand7, Hashim, Yildirim} with VIM, HPM and HFC method respectively. We applied the RBF-DQ method to solve this $Eq.(\ref{Eq5}).$ We construct the residual function as follows:\\
\begin{eqnarray*}
Res(x)=x.\sum_{j=1}^{N} w^{(2)}_{i,j}f_j+2\sum_{j=1}^{N} w^{(1)}_{i,j}f_j+4x(2 e^{f_i}+e^{\frac{f_i}{2}}) ~~~~~i=1,2,\cdots,N.
\end{eqnarray*}
To obtain coefficients $w^{(1)}_{i,j}~and~w^{(2)}_{i,j}$, we put $Res(x)=0$ at the $N$ collocation points, and, then, we construct a set of $N$ equations. 
By solving this set of equations, We have the approximating function $f_i$. 
\begin{table}[htbp]
\centering 
\begin{tabular}{l l l l} 
\hline\hline
$x$ & present Method & Exact value & Error \\
\hline
$0.00$&$0.00000000000$&$0.0000000000$&$0.00e+00$\\
$0.01$&$-0.0001999901$&$-0.0001999900$&$1.24e-10$\\
$0.10$&$-0.0199006604$&$-0.0199006617$&$1.25e-09$\\
$0.50$&$-0.4462871202$&$-0.4462871026$&$1.76e-08$\\
$1.00$&$-1.3862940411$&$-1.3862943611$&$3.20e-07$\\
$2.00$&$-3.2188869094$&$-3.2188758249$&$1.10e-05$\\
$3.00$&$-4.6050645726$&$-4.6051701860$&$1.05e-04$\\
$4.00$&$-5.6664688710$&$-5.6664266881$&$4.21e-05$\\
$5.00$&$-6.5164207344$&$-6.5161930760$&$2.27e-04$\\
$6.00$&$-7.2218661925$&$-7.2218358253$&$3.03e-05$\\
$7.00$&$-7.8240538841$&$-7.8240460109$&$7.87e-06$\\
$8.00$&$-8.3487727119$&$-8.3487745398$&$1.82e-06$\\
$9.00$&$-8..813438533$&$-8.8134384945$&$3.90e-08$\\
$10.00$&$-9.230241034$&$-9.2302410337$&$4.46e-10$\\
\hline
\end{tabular}
\caption{The comparison of the values of $y(x)$ of Lane-Emden equation obtained from present method with $c=1$ and $N=40$ and exact value in Example 5.} 
\label{TabEX5}
\end{table}
\begin{figure}[htbp]
\centering
\includegraphics[width=2in]{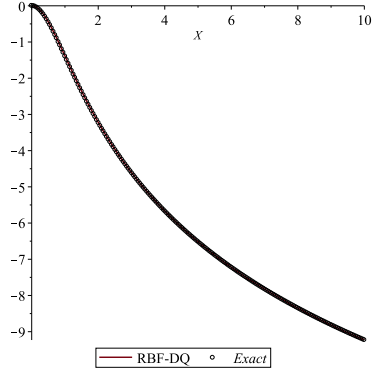}\\
\caption{Graph of equation of example 5 in comparing the presented method and analytic solution }\label{figex5}
\end{figure}
\par
Table \ref{TabEX5} shows the comparison of $y(x)$ obtained by method proposed in this paper and analytic solution Eq.\eqref{Ex5}.\\
In order to compare the present method with the analytic solution,the resulting graph of Eq.\eqref{Ex5} is shown in Fig.\ref{figex5}.
\subsubsection{\textbf{Example 6}}
For $f(x,y)=-6y(x)-4y(x)\ln(y(x))$, $A=1$ and $B=0$, Eq.\eqref{Lane-Emden} has the following form:
\begin{equation}\label{Eq6}
y''(x) + \frac {2}{x} y'(x) -6y(x)-4y(x)\ln(y(x)) = 0, ~~ x>0,
\end{equation}
with conditions:
\begin{eqnarray*}
y(0)=1,\\
y'(0)=0,
\end{eqnarray*}
which has the following analytical solution:
\begin{equation} \label{Ex6}
y(x)= {e^{x^{2}}}.
\end{equation}
This equation has been solved by some researchers; for example \citep{Ramos1, Yildirim} by applying VIM and linearization methods, respectively. If we apply the RBF-DQ method for Eq.\eqref{Eq6}, in this model we have $y(x)\ln(y(x))$ term which increases the order of the calculation; therefore we can use the transformation $y(x)=e^{z(x)}$, in which $z(x)$ is an unknown function for the following equation:
\begin{equation}\label{Ex66}
z''(x) + (z'(x))^{2}+\frac{2}{x} z'(x) -6-4z(x)=0, ~~ x>0,
\end{equation}
with conditions
\begin{eqnarray*}
z(0)=0,\\
z'(0)=0.
\end{eqnarray*}
Now we apply the RBF-DQ method for solving Eq. (\ref{Ex66}) and construct the residual function as follows:
\begin{eqnarray*}
Res(x)= x \sum_{j=1}^{N} w^{(2)}_{i,j}f_j +x (\sum_{j=1}^{N} w^{(1)}_{i,j}f_j)^{2}+x \sum_{j=1}^{N} w^{(1)}_{i,j}f_j -6x-4xf_i=0, ~~~i=1,2,\cdots,N.
\end{eqnarray*}
To obtain coefficients $w^{(1)}_{i,j}~and~w^{(2)}_{i,j}$, we put $Res(x)=0$ at the $N$ collocation points, and, then, we construct a set of $N$ equations.By solving this set of equations, We have the approximating function $f_i$.
\begin{table}[htbp]
\centering 
\begin{tabular}{l l l l l} 
\hline\hline
$x$ & present Method & Exact value & Error \\
\hline
$0.00$&$1.0000000000$&$1.0000000000$&$0.000e+00$\\
$0.01$&$1.0001000050$&$1.0001000050$&$6.05e-22$\\
$0.02$&$1.0004000800$&$1.0004000800$&$2.61e-22$\\
$0.05$&$1.0025031276$&$1.0025031127$&$2.12e-23$\\
$0.10$&$1.0100501670$&$1.0100501671$&$4.75e-25$\\
$0.20$&$1.0408107741$&$1.0408107742$&$5.26e-25$\\
$0.50$&$1.2840254166$&$1.2840254167$&$7.50e-25$\\
$0.70$&$1.6323162199$&$1.6323162200$&$1.07e-24$\\
$0.80$&$1.8964808793$&$1.8964808793$&$1.32e-24$\\
$0.90$&$2.2479079866$&$2.2479079867$&$1.67e-24$\\
$1.00$&$2.7182818284$&$2.7182818285$&$3.009e-19$\\
\hline
\end{tabular}
\caption{The comparison of the values of $y(x)$ of Lane-Emden equation obtained from present method with $c=1$ and $N=35$ and exact value in Example 6.} 
\label{TabEX6}
\end{table}
\begin{figure}
\centering
\includegraphics[width=2in]{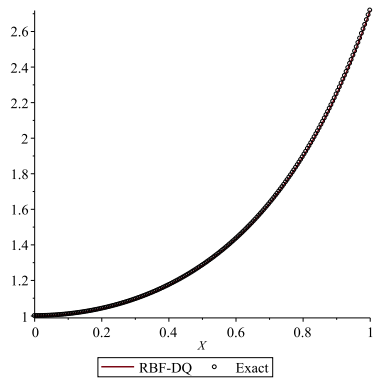}\\
\caption{Graph of equation of example 6 in comparing the presented method and analytic solution }\label{figex6}
\end{figure}
\par
Table\ref{TabEX6} shows the comparison of $y(x)$ obtained by method proposed in this paper and analytic solution, i.e. Eq.\eqref{Eq6}.\\
In order to compare the present method with the analytic solution,the resulting graph of Eq.\eqref{Ex6} is shown in Fig. \ref{figex6}.
\subsubsection{\textbf{Example 7}}
For $f(x,y)=-2(2 x^{2}+3)y$, $A=1$ and $B=0$, Eq.\eqref{Lane-Emden} has the following form:
\begin{equation}\label{Eq7}
y''(x) + \frac {2}{x} y'(x) -2 (2 x^{2}+3)y(x) = 0,~~ x>0,
\end{equation}
with conditions:
\begin{eqnarray*}
y(0)=1,\\
y'(0)=0,
\end{eqnarray*}
which has the following analytical solution:
\begin{equation} \label{Ex7}
y(x)= {e^{x^{2}}}.
\end{equation}
This type of equation has been solved by some researchers, for example \citep{Ramos1, MSH, Yildirim} by using linearization, HPM and VIM methods, respectively. Now we apply RBF-DQ method for Eq.\eqref{Eq7}and construct the residual function as follows:
\begin{eqnarray*}
Res(x)= x \sum_{j=1}^{N} w^{(2)}_{i,j}f_j +2 (\sum_{j=1}^{N} w^{(1)}_{i,j}f_j)^{2}-2x (2x^{2}+3)f_i,~~~~i=1,2,\cdots,N.
\end{eqnarray*}
To obtain coefficients $w^{(1)}_{i,j}~and~w^{(2)}_{i,j}$, we put $Res(x)=0$ at the $N$ collocation points, and, then, we construct a set of $N$ equations. By solving this set, we have the approximating function $f_i$.\\
\begin{table}[htbp]
\centering 
\begin{tabular}{l l l l l} 
\hline\hline
$x$ & present Method & Exact value & Error \\
\hline
$0.00$&$1.0000000000$&$1.0000000000$&$0.00e+00$\\
$0.01$&$1.0001000050$&$1.0001000050$&$9.41e-26$\\
$0.02$&$1.0004000800$&$1.0004000800$&$1.41-25$\\
$0.05$&$1.0025031276$&$1.0025031127$&$2.90e-25$\\
$0.10$&$1.0100501670$&$1.0100501671$&$7.07e-26$\\
$0.20$&$1.0408107741$&$1.0408107742$&$7.09e-25$\\
$0.50$&$1.2840254166$&$1.2840254167$&$1.16e-24$\\
$0.70$&$1.6323162199$&$1.6323162200$&$9.33e-25$\\
$0.80$&$1.8964808793$&$1.8964808793$&$9.64e-25$\\
$0.90$&$2.2479079866$&$2.2479079867$&$4.91e-26$\\
$1.00$&$2.7182818284$&$2.7182818285$&$1.04e-25$\\
\hline
\end{tabular}
\caption{The comparison of the values of $y(x)$ of Lane-Emden equation obtained from present method with $c=1$ and $N=35$ and exact value in Example 7.} 
\label{TabEX7}
\end{table}
\begin{figure}[htbp]
\centering
\includegraphics[width=2in]{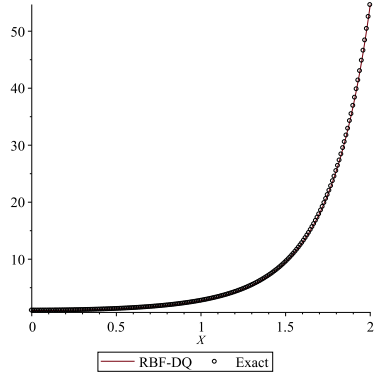}\\
\caption{Graph of equation of example 7 in comparing the presented method and analytic solution }\label{figex7}
\end{figure}
\par
Table\ref{TabEX7} shows the comparison of $y(x)$ obtained by method proposed in this paper and analytic solution, i.e. Eq.\eqref{Eq7}.\\
In order to compare the present method with the analytic solution,the resulting graph of Eq.\eqref{Ex7} is shown in Fig.\ref{figex7}.
\subsubsection{\textbf{Example 8}}
For $f(x,y)=xy(x),~ h(x)=x^5-x^4+44x^2-30x,~ A=0$ and $B=0$, Eq.\eqref{Lane-Emden} will be one of the Lane-Emden type equations that is absorbing to solve
\begin{equation}\label{Eq8}
y''(x)+\frac{8}{x}y'(x)+xy(x)=x^5-x^4+44x^2-30x, x\geq 0,
\end{equation}
subject to the boundary conditions
\begin{eqnarray*}
y(0)=0,\\
y'(0)=0,
\end{eqnarray*}
which has the following analytical solution:
\begin{equation}\label{Ex8}
y(x)=x^4+x^3.
\end{equation}
This type of equation has been solved by \citep{Parand7, Ramos1, Bataineh, Ozis1} with HFC, linearization, HAM and HPM method respectively. We applied the RBF-DQ method to solve this Eq.\eqref{Eq8}. We construct the residual function as follows:\\
\begin{eqnarray*}
Res(x)= x \sum_{j=1}^{N} w^{(2)}_{i,j}u_j +8\sum_{j=1}^{N} w^{(1)}_{i,j}u_j+x (xu_i-x^5+x^4-44x^2+30x), ~~~~~i=1,2,\cdots,N.
\end{eqnarray*}
To obtain coefficients $w^{(1)}_{i,j}~and~w^{(2)}_{i,j}$, we put $Res(x)=0$ at the $N$ collocation points, and, then, we construct a set of $N$ equations. By solving this set of equations, We have the approximating function $f_i$.
\begin{table}[htbp]
\centering 
\begin{tabular}{l l l l} 
\hline\hline
$x$ & present Method & Exact value & Error \\
\hline
$0.00$ &$0.0000000000$ &$0.0000000000$ &$0.00e+00$\\
$0.01$ &$-9.899763744$ &$-0.000000900$ &$2.36e-11$\\
$0.10$ &$-0.000899999$ &$-0.000900000$ &$3.68e-10$\\
$0.50$ &$-0.062500014$ &$-0.062500000$ &$1.44e-08$\\
$1.00$ &$0.0000003209$ &$0.0000000000$ &$3.20e-07$\\
$2.00$ &$8.0000269820$ &$8.0000000000$ &$2.69e-06$\\
$3.00$ &$53.999201677$ &$54.000000000$ &$7.98e-04$\\
$4.00$ &$192.00521430$ &$192.00000000$ &$5.21e-03$\\
$5.00$ &$499.99587333$ &$500.00000000$ &$4.12e-03$\\
$6.00$ &$1079.9933014$ &$1080.0000000$ &$6.69e-03$\\
$7.00$ &$2058.0030011$ &$2058.0000000$ &$3.001e-03$\\
$8.00$ &$3583.9998744$ &$3584.0000000$ &$1.25e-04$\\
$9.00$ &$5831.9999966$ &$5832.0000000$ &$3.33e-06$\\
$10.0$ &$9000.0000000$ &$9000.0000000$ &$4.04e-09$\\
\hline
\end{tabular}
\caption{The comparison of the values of $y(x)$ of Lane-Emden equation obtained from present method with $c=1$ and $N=45$ and exact value in Example 8.} 
\label{TabEX8}
\end{table}
\begin{figure}[htbp]
\includegraphics[width=3in]{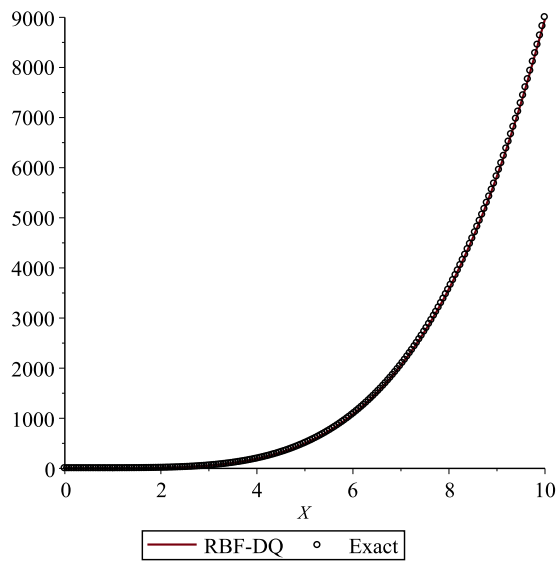}\\
\caption{Graph of equation of example 8 in comparing the presented method and analytic solution }\label{figex8}
\end{figure}
\par
Table \ref{TabEX8} shows the comparison of $y(x)$ obtained by method proposed in this paper and analytic solution, i.e. Eq.\eqref{Eq8}.\\
In order to compare the present method with the analytic solution,the resulting graph of Eq.\eqref{Ex8} is shown in Fig. \ref{figex8}.\\
\subsubsection{\textbf{Example 9}} 
For $f(x,y)=y(x),~ h(x)=6+12x+x^2+x^3,~ A=0$ and $B=0$, Eq.\eqref{Lane-Emden} will be one of the Lane-Emden type equations that is absorbing to solve
\begin{equation}\label{Eq9}
y''(x)+\frac{2}{x}y'(x)+y(x)=6+12x+x^2+x^3, x\geq 0,
\end{equation}
subject to the boundary conditions
\begin{eqnarray*}
y(0)=0,\\
y'(0)=0,
\end{eqnarray*}
which has the following analytical solution:
\begin{equation}\label{Ex9}
y(x)=x^4-x^3.
\end{equation}
This equation has been solved by \citep{Parand7, Ramos1, Bataineh, Ozis1} with HFC, linearization, HAM and HPM method respectively. We applied the RBF-DQ method to solve equation Eq.\eqref{Eq9}.Therefore,we construct the residual function as follows:
\begin{eqnarray*}
Res(x)= x \sum_{j=1}^{N} w^{(2)}_{i,j}u_j +8\sum_{j=1}^{N} w^{(1)}_{i,j}u_j+x (xu_i-x^5+x^4-44x^2+30x), ~~~~~i=1,2,\cdots,N.
\end{eqnarray*}
To obtain coefficients $w^{(1)}_{i,j}~and~w^{(2)}_{i,j}$, we put $Res(x)=0$ at the $N$ collocation points, and, then, we construct a set of $N$ equations.By solving this set of equations, We have the approximating function $u_i$.
\begin{table}[htbp]
\centering 
\begin{tabular}{l l l l} 
\hline\hline
$x$ & present Method & Exact value & Error \\
\hline
$0.00$ &$0.0000000000$ &$0.0000000000$ &$0.00e+00$\\
$0.01$ &$0.0001010000$ &$0.0001010000$ &$8.85e-12$\\
$0.10$ &$0.0110000000$ &$0.0110000000$ &$1.36e-11$\\
$0.50$ &$0.3749999976$ &$0.3750000000$ &$2.39e-09$\\
$1.00$ &$2.0000000323$ &$2.0000000000$ &$3.23e-08$\\
$2.00$ &$12.000002499$ &$12.000000000$ &$2.49e-06$\\
$3.00$ &$35.999914812$ &$36.000000000$ &$8.51e-05$\\
$4.00$ &$80.000486174$ &$80.000000000$ &$4.86e-05$\\
$5.00$ &$149.99982164$ &$150.00000000$ &$1.78e-04$\\
$6.00$ &$251.99924878$ &$252.00000000$ &$7.51e-04$\\
$7.00$ &$392.00029168$ &$392.00000000$ &$2.91e-04$\\
$8.00$ &$575.99998728$ &$576.00000000$ &$1.27e-05$\\
$9.00$ &$809.99999967$ &$810.00000000$ &$3.26e-07$\\
$10.0$ &$1100.0000000$ &$1100.0000000$ &$1.52e-10$\\
\hline
\end{tabular}
\caption{The comparison of the values of $y(x)$ of Lane-Emden equation obtained from present method with $c=1$ and $N=45$ and exact value in Example 9.} 
\label{TabEX9}
\end{table}
\begin{figure}[htbp]
\centering
\includegraphics[width=2.2in]{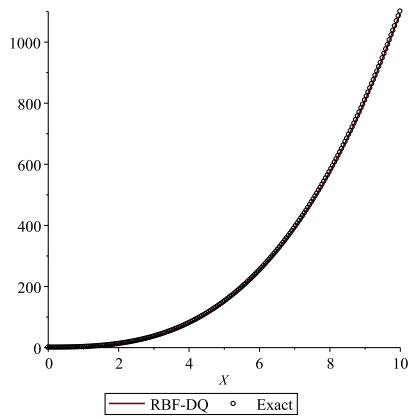}\\
\caption{Graph of y(x) of example 9 in comparing the presented method and analytic solution }\label{figex9}
\end{figure}
\par
Table\ref{TabEX9} shows the comparison of $y(x)$ obtained by method proposed in this paper and analytic solution, i.e. Eq.\eqref{Eq9}.\\
In order to compare the present method with the analytic solution,the resulting graph of Eq.\eqref{Ex9} is shown in Fig. \ref{figex9}.\\
\section{ Conclusion}
The Lane–Emden equation describes a variety of phenomena in theoretical physics and astrophysics, including aspects of stellar structure, the thermal history of a spherical cloud of gas, isothermal gas spheres, and thermionic currents \citep{Chandra}.
Since then the equation has been a center of attraction. The main problem in this direction is the accuracy and range of applicability of these approaches.\par
The main goal of this paper was to introduce a RBF-DQ method to construct an approximation to the solution of non-linear Lane–Emden equation in a semi-infinite interval $[0,\infty)$, Which is meshfree characteristic and does not require mesh generation. Further, RBF-DQ, in contrast to classical DQ, can be directly applied to irregular domains.\par
Our results have better accuracy with compared to other results. The validity of the method is based on the assumption that it converges when increasing the number of collocation points. A comparison is made among the exact solution and the numerical solution of Horedt\citep{Horedt} and series solutions of Wazwaz\citep{Wazwaz1}, Parand et al.\citep{Parand8} and the current work. The results of the present method for this type of problem clearly indicate that our method is accurate.\\

\end{document}